\documentclass[a4paper,11pt,twoside]{article}
\usepackage{pst-fill,pst-grad}
\usepackage{textcomp} 
\usepackage[english]{babel}  
\usepackage[utf8]{inputenc}    
\usepackage[titletoc]{appendix}
\usepackage{titlesec}  
\usepackage{graphicx} 
\usepackage{amsmath}  
\usepackage{float} 
\usepackage{fancyhdr}  
\usepackage[matrix,arrow,curve]{xy} 
\usepackage{pstricks} 
\usepackage{amsmath,amsfonts,verbatim,afterpage,theorem,euscript,mathrsfs,amssymb}
\usepackage{amsfonts}
\usepackage{amssymb}
\usepackage{array} 
\usepackage{dsfont}
\usepackage[colorlinks=true,linkcolor=blue,citecolor=red]{hyperref}
\usepackage{authblk}
\usepackage{color}
\newtheorem{Proposition}{Proposition}[section]

\newtheorem{Theoreme}{Theorem}[section]

\newtheorem{Corollaire}{Corollary}[section]
\newtheorem{Remarque}{Remark}

\def \vu{\vec{u}}
\def \vb{\vec{b}}

\def \P{\mathbb{P}}

\def \R{\mathbb{R}}
\def \Rn{\mathbb{R}^n}
\def \Rt{\mathbb{R}^3}

\def \finpv{\hfill $\blacksquare$} 
\def \pv{{\bf{Proof.}}~} 

\def \ds{\displaystyle}

\title{\bf From non-local to local Navier-Stokes equations}

\author[1]{ Oscar Jarr\'in\footnote{corresponding author: oscar.jarrin@udla.edu.ec}}
\author[2]{ Geremy Loacham\'in\footnote{geremy.loachamin@uni.lu}}

\affil[1]{\scriptsize Escuela de Ciencias Físicas y Matemáticas, Universidad de Las Américas, Vía a Nayón, C.P.170124, Quito, Ecuador.} 
\affil[2]{\scriptsize Faculty of Science, University of Luxembourg, Maison du Nombre, 6 Avenue de la Fonte, Esch-sur-Alzette L-4364, Luxembourg.}
\date{\today}
\begin{document}
	\maketitle	 
	%%%%%%%%%%%%%%%%%%%%%%%%%%%%%%%%%%%%%%%%%%%%%%
\begin{abstract}
	Inspired by some experimental (numerical) works on fractional diffusion PDEs, we develop a rigorous framework to prove that solutions to the fractional Navier-Stokes equations, which involve the fractional Laplacian operator $(-\Delta)^{\frac{\alpha}{2}}$ with $\alpha<2$, converge to a solution of the classical case, with $-\Delta$, when $\alpha$ goes to $2$. Precisely, in the setting of mild solutions, we prove uniform convergence in both the time and spatial variables and derive a precise convergence rate, revealing some  phenomenological effects. Finally, our results are also generalized to the coupled setting of the Magnetic-hydrodynamic  system. \\[3mm]
\end{abstract}
	%\tableofcontents
	%%%%%%%%%%%%%%%%%%%%%%%%%%%%%%%%%%%%%%%%%%%%%%
	%%%%%%%%%%%%%%%%%%%%%%%%%%%%%%%%%%%%%%%%%%%%%%

\section{Introduction} 
For a velocity field $\vu: [0,+\infty)\times \Rt \to \Rt
$, and for a pressure term $p:[0,+\infty)\times \Rt \to \R$,  we deal with  the three-dimensional, incompressible and generalized Navier-Stokes equations  in the whole space $\Rt$:
\begin{equation}\label{NS-Intro}
	\partial_t \vu = -\nu (-\Delta)^{\alpha/2} \vu - (\vu \cdot \vec{\nabla})\vu -\vec{\nabla} p, \quad \text{div}(\vu)=0, \qquad  1<\alpha \leq 2.
\end{equation}

When $1<\alpha<2$, the diffusion term is given by the fractional Laplacian operator, which is easily defined in the Fourier variable by the symbol $|\xi|^\alpha$. Moreover, in the spatial variable, we have 
\[ (-\Delta)^{\alpha/2} \vu(t,x)= C_{\alpha}\, {\bf p.v.} \int_{\Rt} \frac{\vu(t,x)- \vu(t,y)}{|x-y|^{3+\alpha}} dy,  \]
where $C_\alpha>0$ is a constant depending on $\alpha$, and  ${\bf p.v.}$ denotes the principal value.  The non-local behavior of this operator allows us to call the equations (\ref{NS-Intro}) the non-local Navier-Stokes equations. By contrast, when $\alpha=2$, the diffusion term is given by the classical Laplacian operator, and we shall refer to the classical (or local) Navier-Stokes equations.  With a minor loss of generality, we shall set the viscosity constant $\nu$ equal to one. 

\medskip

Numerical solutions to the classical Navier-Stokes equations (when $\alpha=2$) for engineering problems, turbulent fluid flows, and geophysical phenomena are not completely possible at present, see for example \cite{Dubois,Pope}. In addition, the mathematical theory of global existence and regularity of solutions to these equations remains one of the most challenging open questions in mathematical analysis, as discussed in \cite{PGLR,Temam}. In this context, the fractional Navier-Stokes equations (when $1<\alpha<2$) have been employed as a relevant modification of the classical equations to gain a better understanding of these mathematical and computational difficulties \cite{Holst,Meerschaert,Nan,Olson-Titi}. 

\medskip

In equation (\ref{NS-Intro}), for each  $1<\alpha<2$ fixed, we have an associated fractional Navier-Stokes equation, and we will denote its solution by  $(\vu_{\alpha}, p_\alpha)$. The main objective of this note is to sharply study the \emph{dynamics} of the family of solutions $(\vu_{\alpha}, p_\alpha)_{1<\alpha<2}$ when the parameter $\alpha$ tends to $2$.  

\medskip

This question is not only interesting from the theoretical point of view, but has also been pointed out in some experimental works involving fractional Burgers equations \cite{Funaki} and a fractional transport-type equation \cite{Zaslavshy}. More precisely, these  \emph{numerical} studies show that solutions to fractional equations behave as solutions to the classical ones (involving the Laplace operator) when $\alpha$ is sufficiently close to $2$. Inspired by these previous works, we aim to develop a \emph{rigorous framework} to study the convergence
\begin{equation}\label{Convergence-Intro}
	(\vu_{\alpha}, p_\alpha) \to (\vu_2, p_2), \quad \mbox{when} \quad \alpha \to 2^{-}, 
\end{equation}
where $(\vu_2, p_2)$ denotes a solution to the classical Navier-Stokes equations.

\medskip

It is also worth mentioning that this question has been studied for some \emph{elliptic} equations, including the nonlinear Schr\"odinger equation \cite{Bieganowski} and the fractional $p$-Laplacian problem \cite{Fernandez-Salort}. In these works, the authors mainly used variational methods and concentration-compactness principles to prove the convergence of \emph{weak solutions} of the fractional problem to the classical problem.  Specifically, in \cite{Bieganowski}, this convergence was proven in the strong topology of the space $L^{2}_{loc}(\R^n)$ (with $n\geq 3$), whereas in \cite{Fernandez-Salort}, the authors used the (more technical) notion of $\Gamma$-convergence. 

\medskip

For the \emph{parabolic} setting of equation (\ref{NS-Intro}), a first interesting study of the   convergence (\ref{Convergence-Intro}) was done in \cite{Dlotko} for both the  two-dimensional case and the  three-dimensional case, in the setting of a bounded and smooth domain $\Omega \subset \R^n$ (with $n=2,3$), and for the sub-critical  case when $\alpha >2$.  Specifically, in the two-dimensional case,  it is proven that a family of  solutions  $(\vu_\alpha)_{2<\alpha \leq 5/2}$ to equations (\ref{NS-Intro}) convergences (when $\alpha \to 2^{+}$) in the \emph{weak topology} of the space $L^2(\Omega)$ to a \emph{weak Leray's solution}  $\vu_2$ of the classical Navier-Stokes equations.  This result does not   fulfill for the three-dimensional case, due to the  lost of some  key tools only available in $2D$ to handle the nonlinear term. Thus, in this case it is proven that solutions to the regularized equations
\[ \partial_t  \vu  = \big(\Delta - \varepsilon (-\Delta)^{\alpha/2}\big) \vu - (\vu \cdot \vec{\nabla})\vu -\vec{\nabla} p, \quad \text{div}(\vu)=0, \qquad \alpha>5/2,\]
converge (when $\varepsilon \to 0^{+}$) to a \emph{weak Leray's solution} $\vu_2$ of the classical Navier-Stokes equations in the \emph{weak topology} of the energy space. The ideas to proof these results are mainly based on sharp \emph{a priori} energy estimates, the weak formulation of solutions and concentration-compactness arguments.  

\medskip

In this paper, we will employ a completely different approach (following some of the ideas presented in our previous work \cite{JarrinLoachamin}), which is  principally based on the explicit structure of \emph{mild solutions}. This approach and some sharp computations in the Fourier level (since we consider here the whole space $\Rt$) allow us to prove a \emph{uniform convergence} (\ref{Convergence-Intro}) in the \emph{strong topology} of the $L^{\infty}_{tx}$-space, and in the super-critical case when $\alpha <2$.  Of course, under minor technical adaptations, our results hold for the sub-critical case when $\alpha>2$. Moreover, our approach 
 also allows us to  derive an \emph{explicit convergence rate} which highlights some interesting phenomena, that we shall expose in detail below. 
\medskip 
 
{\bf The main result.} We focus on the initial value problem for both the non-local ($1<\alpha<2$) and local ($\alpha=2$) Navier-Stokes equations:
\begin{equation}\label{NS} 
	\begin{cases}\vspace{2mm}
		\partial_t \vu_\alpha =- (-\Delta)^{\alpha/2} \vu_\alpha - (\vu_\alpha \cdot \vec{\nabla})\vu_\alpha -\vec{\nabla} p_\alpha, \quad  \text{div}(\vu_\alpha)=0, \qquad  1<\alpha \leq 2,  \\ 
		\vu_\alpha(0,\cdot)=\vu_{0,\alpha},  
	\end{cases}
\end{equation}
where $\vu_{0,\alpha}:\Rt \to \Rt$ denotes the (divergence-free) initial datum. Recall that \emph{mild solutions} to equations (\ref{NS}) are obtained using Banach's contraction principle by solving the following integral equation (due to Duhamel's formula)
\begin{equation}\label{Mild}
	\vu_\alpha(t,\cdot)=e^{-(-\Delta)^{\alpha/2} t} \, \vu_{0,\alpha} - \int_{0}^{t} e^{-(-\Delta)^{\alpha/2} (t-\tau)} \, \P\left((\vu_\alpha \cdot \vec{\nabla})\vu_\alpha\right)(\tau, \cdot) d \tau, \qquad 1<\alpha \leq 2.
\end{equation}

 In this expression, for $1<\alpha<2$, we have $e^{-(-\Delta)^{\alpha/2} t} f= h_\alpha(t,\cdot)\ast f$, where the kernel $h_\alpha(t,x)$ is the fundamental solution to the fractional heat equation $\partial_t h_\alpha + (-\Delta)^{\frac{\alpha}{2}}  h_\alpha=0$ when $t>0$, and $h_\alpha(0,\cdot)=\delta_0$ where $\delta_0$ is the Dirac mass at the origin. For $\alpha=2$, we have $e^{\Delta t} f= h(t,\cdot)\ast f$, where $h(t,x)$ is the well-known heat kernel.
 
 \medskip
 
The operator $\P$ stands for Leray's projector, and it is well-known that the pressure $p_\alpha$ can be easily deduced from the velocity $\vu_{\alpha}=(u_{\alpha,1},u_{\alpha,2},u_{\alpha,3})$ due to the divergence-free property of this latter. Then, we have 
\begin{equation}\label{Pressure}
	p_\alpha= \frac{1}{-\Delta} \text{div}\left((\vu_\alpha \cdot \vec{\nabla})\vu_\alpha\right)= \sum_{i,j=1}^{3} \mathcal{R}_{i}\mathcal{R}_{j} (u_{\alpha,i} \, u_{\alpha,j}),
\end{equation}
where $\mathcal{R}_{i}=\frac{\partial_i}{\sqrt{-\Delta}}$ denotes the Riesz transform. 
 
 \medskip
 
In the setting of non-homogeneous Sobolev spaces $H^s(\Rt)$, with $s>1/2$, the local well-posedness theory for mild solutions to the classical Navier-Stokes equations is a well-known issue \cite{Chemin}. In our next proposition, for the sake of completeness of this article, we revisit this result for the generalized case of equation (\ref{NS}) in the space $H^s(\Rt)$  with $s>3/2$. The (technical) constraint $s>3/2$ will be  useful later  to prove our \emph{key tool}, given in Proposition \ref{Key-Lemma} below, for the study of the convergence (\ref{Convergence-Intro}). Precisely, our  approach is  based on a sharp study of the  convergence  $h_\alpha(t,\cdot) \to h(t,\cdot)$ when $\alpha \to 2^{-}$, for the kernels in the mild formulation (\ref{Mild}).  

\medskip 

We emphasize that the proof of the proposition below is classical, but we also aim to determine how the existence time of the mild solution $\vu_{\alpha}$, denoted by $T_\alpha$, explicitly depends on the parameter $\alpha$. 
\begin{Proposition}\label{Prop:LWP} Let $1<\alpha \leq 2$ be fixed. Let $s>3/2$, and let $\vu_{0,\alpha} \in H^s(\Rt)$ be a divergence-free initial datum. There exists a time 
\begin{equation}\label{Talpha}
0<T_\alpha = \frac{1}{2} \left( \frac{1- \frac{1}{\alpha}}{4C \| \vu_{0,\alpha} \|_{H^s}} \right)^{\frac{\alpha}{\alpha-1}}, 
\end{equation}	
where $C>0$ is a generic constant, and there exists a unique mild solution $\vu_\alpha$ to equation (\ref{NS}) such that 
\[ \vu_\alpha \in \mathcal{C}\big([0,T_\alpha], H^s(\Rt)\big) \quad \mbox{and} \quad  p_\alpha \in \mathcal{C}\big([0,T_\alpha], H^s(\Rt)\big). \]
\end{Proposition}	
\begin{Remarque} Note that $0<T_\alpha$ as long as $1<\alpha$.
\end{Remarque}	

Once we have stated this proposition, we rigorously studied the convergence presented in (\ref{Convergence-Intro}). For the non-local  case (when $1<\alpha<2$), we consider a family of initial data $(\vu_{0,\alpha})_{1<\alpha <2} \subset H^s(\Rt)$ and the resulting family of solutions $(\vu_{\alpha},p_\alpha)_{1<\alpha<2}$  obtained in Proposition \ref{Prop:LWP}. Similarly, for the local case (when $\alpha=2$), we consider the initial datum $\vu_{0,2}\in H^s(\Rt)$ and its associated solution $(\vu_2,p_2)\in \mathcal{C}([0,T_2], H^s(\Rt))$.

\medskip

Our starting point is to assume the following convergence of initial data in the strong topology of the space $H^s(\Rt)$: 
\begin{equation}\label{Conv-Data}
	\vu_{0,\alpha} \to \vu_{0,2}, \quad \alpha \to 2^{-}. 
\end{equation}

This assumption yields the following important facts. On one hand, this convergence will allow us to find a quantity $0<\varepsilon<\ll 1$ and a time $T_0$, depending only on $\varepsilon$, such that 
\begin{equation}\label{Lower-bound-time}
	T_0\leq T_\alpha, \quad \mbox{for all} \quad 1+\varepsilon < \alpha \leq 2. 
\end{equation}

See Appendix \ref{AppendixA} for a rigorous justification of this fact. Consequently, for $1+\varepsilon<\alpha \leq 2$, each solution $(\vu_{\alpha},p_\alpha)$ is defined at least on the time interval $[0,T_0]$, and this fact will be used when studying (\ref{Convergence-Intro}). 

\medskip

On the other hand, since $s>3/2$, the space $H^s(\Rt)$ is continuously embedded in the space $L^\infty(\Rt)$, and the convergence (\ref{Conv-Data}) also holds in $L^\infty(\Rt)$. Thus, 
for the family of velocities $\vu_{\alpha}$, we shall prove the  following \emph{uniform convergence}:
\begin{equation}\label{Convergence-velocities}
	\vu_{\alpha}\to \vu_2, \quad \alpha \to 2^{-}, \quad \mbox{in}\quad L^{\infty}([0,T_0]\times \Rt).
\end{equation}

Recall that the pressures $p_\alpha$ are defined through Riesz transforms and the velocities $\vu_{\alpha}$ in the expression (\ref{Pressure}). Nevertheless, since Riesz transforms are not bounded in the $L^{\infty}-$space, we need to consider the larger space $BMO(\Rt)$. See \cite[Chapter $3$]{Grafakos} for a definition and some properties of this space. In this setting, convergence (\ref{Convergence-velocities}) yields to prove:
\begin{equation}\label{Convergence-pressures}
	p_{\alpha}\to p_2, \quad \alpha \to 2^{-}, \quad \mbox{in}\quad L^{\infty}([0,T_0], BMO(\Rt)).
\end{equation}

Furthermore, our main contribution is to quantify how fast the convergences (\ref{Convergence-velocities}) and (\ref{Convergence-pressures}) hold.  To this end, for a  parameter  $\kappa>0$ fixed,  we shall assume a convergence rate of initial data (see (\ref{Conv-Rate-Data}) below) which is measured in terms of  $\kappa$. We thus aim to study when this convergence rate persists for solutions. In this context, we present our main result:
\begin{Theoreme}\label{Th-Main} Let $\ds{(\vu_{0,\alpha})_{1+\varepsilon<\alpha \leq 2}}$ be an initial data family, where $\vu_{0,\alpha} \in H^s(\Rt)$ with $s>3/2$. Then,  let 
\[  (\vu_\alpha, p_\alpha)_{1+\varepsilon<\alpha\leq 2} \subset  \mathcal{C}\big([0,T_0], H^s(\Rt)\big), \]	
be the associated family of solutions to  equation (\ref{NS}), obtained in Proposition \ref{Prop:LWP}. 
	
\medskip 

We assume the convergence given in (\ref{Conv-Data}). Moreover, for a parameter $\kappa>0$, we assume the convergence rate of initial data 
\begin{equation}\label{Conv-Rate-Data}
	\|  \vu_{0,\alpha} - \vu_{0,2} \|_{L^\infty} \leq {\bf c} (2-\alpha)^\kappa,
\end{equation}
where ${\bf c}>0$ is a generic constant. There exists a constant $0<{\bf C}(T_0) \sim 1+T_0 +T^{2}_{0}$, depending on the initial datum $\vu_{0,2}$, the quantity $\varepsilon$, the constant ${\bf c}$ and the time $T_0$, such that for  all  $1+\varepsilon < \alpha < 2$  the following estimate holds:
\begin{equation}\label{Conv-Rate-Sol}
 		\sup_{0\leq t \leq T_0}  \Big( \| \vu_\alpha(t,\cdot) - \vu_2(t,\cdot)\|_{L^\infty} + \| p_\alpha(t,\cdot) - p_2(t,\cdot)\|_{BMO}\Big) \leq   {\bf C}(T_0)  \Big( (2-\alpha)+(2-\alpha)^\kappa \Big). 
\end{equation}
\end{Theoreme}
 
Some remarks have been provided in order here. The uniform convergence (\ref{Conv-Rate-Sol}) is stronger than the ones obtained in the aforementioned works \cite{Bieganowski,Dlotko,Fernandez-Salort}. Moreover, in contrast to these works, we also derive a convergence rate (when $\alpha \to 2^{-}$) of the order $(2-\alpha)+(2-\alpha)^\kappa$. 
\begin{Remarque} In this last expression, we observe that increasing values of $\kappa$ makes our hypothesis (\ref{Conv-Rate-Data}) strong, but \emph{this fact does not persist for solutions}.  Precisely, for $\alpha$ sufficiently close to $2$,  we obtain
\[ (2-\alpha)+(2-\alpha)^\kappa \sim (2-\alpha), \quad \mbox{when} \quad 1<\kappa, \]
and consequently, the convergence rate of solutions is \emph{slower} than the one of initial data.
\end{Remarque}	

To understand this \emph{unexpected  phenomenon}, let us recall that mild solution to equation (\ref{NS}) are given in expression (\ref{Mild}), where the main difference between the fractional case and the classical one are the kernels $h_\alpha(t,\cdot)$ and $h(t,\cdot)$ respectively. In Proposition \ref{Key-Lemma} below, we rigorously prove the convergence $h_\alpha(t,\cdot)\to h(t,\cdot)$ (when $\alpha \to 2^{-}$) with a \emph{optimal} convergence rate of the order $(2-\alpha)$.  Therefore, the convergence rate of solutions is given by a \emph{competition} between the \emph{assumed} convergence rate of initial date and the \emph{phenomenological} convergence rate of the kernels in the mild formulation.  

\begin{Remarque}
	 In the particular case of  the same initial data for whole the family of equations (\ref{NS}): $\vu_{0,\alpha}=\vu_{0,2}$ for all $1+\varepsilon< \alpha <2$, the estimate (\ref{Conv-Rate-Sol}) becomes 
	\[  \sup_{0\leq t \leq T_0}  \big( \| \vu_\alpha(t,\cdot) - \vu_2(t,\cdot)\|_{L^\infty} + \| p_\alpha(t,\cdot) - p_2(t,\cdot)\|_{BMO}\big)  \leq {\bf C}(T_0) (2-\alpha), \]
	where the convergence rate is purely determined by the convergence of the kernels $h_\alpha(t,x)\to h(t,x)$.
\end{Remarque}	
In the case of small initial data,  it is well known that mild solutions to equation (\ref{NS}) are global in time, see \cite[Theorem $7.3$]{PGLR}. In this setting, we have
\begin{Corollaire}\label{Corollary:global} Under the same hypothesis as in Theorem \ref{Th-Main}, assume that 
\begin{equation}\label{Small-data}
\sup_{1+\varepsilon<\alpha \leq 2}\| \vu_{0,\alpha} \|_{H^s} \ll 1.
\end{equation}	
Then, for all $1+\varepsilon < \alpha < 2$, our main  estimate (\ref{Conv-Rate-Sol}) writes down as 
\begin{equation}\label{Conv-Rate-Sol-global}
\sup_{0\leq t \leq T} \Big( \| \vu_\alpha(t,\cdot) - \vu_2(t,\cdot) \|_{L^\infty} + \| p_\alpha(t,\cdot) - p_2(t,\cdot) \|_{BMO}\Big) \\
\leq \,  {\bf C}(T) \Big( (2-\alpha)+ (2-\alpha)^\kappa\Big), 
\end{equation}
where,  for all time $0<T<+\infty$ we have ${\bf C}(T) \sim 1+T+T^2$. Moreover, the limit $\vu_2$ is also a Leray's solution to the classical Navier-Stokes equations. 
\end{Corollaire}	

The convergence result presented in Theorem \ref{Th-Main} also allows us to study the convergence from the non-local to the local Navier-Stokes equation in the space $L^{p}((0,T_0), L^q(\Rt))$. In the above estimate, it is interesting to observe that the convergence rate depends on  the parameters $\kappa$ and $q$ but not on the parameter $p$. 
	\begin{Corollaire}\label{Corollary:LpLq} Under the same hypothesis as in Theorem \ref{Th-Main}, for  $1\leq p \leq +\infty$ and $2<q<+\infty$ the estimate holds: 
	\begin{equation}\label{Conv-Rate-Sol-Lp}
			\left\| \vu_\alpha - \vu_2 \right\|_{L^{p}_{t}L^{q}_{x}}  + \left\| p_\alpha - p_2 \right\|_{L^{p}_{t}L^{q}_{x}} \\
			\leq \,  {\bf C}_{p,q}(T_0) \Big((2-\alpha) + (2-\alpha)^\kappa \Big)^{1-1/q},
	\end{equation}
	with ${\bf C}_{p,q}(T_0) \sim (1+T_0+T^{2}_{0})$. 
	\end{Corollaire}	

  Finally, it is worth mentioning that Theorem \ref{Th-Main} also holds for  to the two-dimensional case, where the regularity constraint $s>3/2$ is relaxed to $s>1$. In this sense, we complete the previous work \cite{Dlotko} with a non-local to local convergence result for $2D$ Navier-Stokes mild solutions. 
  
  \medskip
  
  On the other hand, Theorem \ref{Th-Main} can be generalized to some relevant coupled systems in fluid dynamics, such as the Magneto-hydrodynamic (MHD) equations. See Appendix \ref{AppendixB} for all details.

 \medskip

To conclude this section, we would like to make some final comments: as mentioned, the strategy developed to prove Theorem \ref{Th-Main} is strongly based on mild solutions of the equations (\ref{NS}). In future research, we aim to develop a different approach to study the convergence (\ref{Convergence-Intro}) in the setting of Leray's solutions. Moreover, by following some of the ideas in \cite{Bieganowski, Fernandez-Salort}, we think it would be interesting to study this convergence in the elliptic case of \emph{stationary} (time-independent) solutions. 

\medskip 

{\bf Organization of the article.}  Section \ref{Sec:Prelimonaries} is essentially devoted to the proof of the key Proposition \ref{Key-Lemma}. In Section \ref{Sec:LWP}, for the sake of completeness, we provide a brief proof of Proposition \ref{Prop:LWP}. Finally, in Section \ref{Sec:NL-Ll}, we prove our main results: Theorem \ref{Th-Main} and its Corollary \ref{Corollary:LpLq}.

\section{Preliminaries: Non-local to local heat equation}\label{Sec:Prelimonaries} 
Recall that the fractional kernel $h_\alpha(t,x)$ is easily defined in the Fourier level by the expression $\ds{\widehat{h_\alpha}(t,\xi)= e^{-t\, \vert \xi \vert^\alpha}}$, while in the classical case, we have  $\ds{\widehat{h}(t,\xi)= e^{-t\, \vert \xi \vert^2}}$. 

\medskip

In the following proposition, we study the strong convergence of the kernel $h_\alpha(t,x)$ to the heat kernel $h(t,x)$, when $\alpha \to 2^{-}$. This result will be our key tool in the sequel. 
\begin{Proposition}\label{Key-Lemma} Let $s>3/2$. There exists a constant $C=C_s>0$  such that, for all $1<\alpha < 2$ and for all time $0<T<+\infty$, the following estimate from above holds: 
	\begin{equation}\label{Above}
	\sup_{0\leq t \leq T} \left\Vert h_\alpha(t,\cdot)-h(t,\cdot) \right\Vert_{H^{-s}}\leq  C \, T ( 2- \alpha ).
	\end{equation}
Moreover, there exists a constant $c=c_s \leq C$, and there exists  quantity $0<\varepsilon_1 \ll 1$ such that for all $1+\varepsilon_1 < \alpha < 2$ the estimate from below holds:
\begin{equation}\label{Below}
c \, \frac{T}{2} (2-\alpha) \leq \sup_{0\leq t \leq T} \left\Vert h_\alpha(t,\cdot)-h(t,\cdot) \right\Vert_{H^{-s}}.
\end{equation}
\end{Proposition}
\pv We begin by verifying that the expression $\ds{ \left\Vert h_\alpha(t,\cdot)-h(t,\cdot) \right\Vert^{2}_{H^{-s}}}$ is a continuous function of $t$. For $0\leq t_0, t\leq T$, we have 
\begin{equation*} 
	\left\Vert h_\alpha(t,\cdot)-h(t,\cdot) \right\Vert^{2}_{H^{-s}} - \left\Vert h_\alpha(t_0,\cdot)-h(t_0,\cdot) \right\Vert^{2}_{H^{-s}}  
	=  \int_{\Rt} \left( \left\vert  e^{-\vert \xi \vert^{\alpha} t} -e^{-\vert \xi \vert^2 t}\right\vert^2- \left\vert  e^{-\vert \xi \vert^{\alpha} t_0} -e^{-\vert \xi \vert^2 t_0}\right\vert^2 \right)\frac{d \xi}{(1+\vert \xi \vert^2)^{s}}.
\end{equation*}

As  $s>3/2$, we have $\ds{\int_{\Rt}\frac{d \xi}{(1+\vert \xi \vert^2)^{s}} <+\infty}$, and  a direct application of the well-known dominated  convergence theorem yields
\[ \lim_{t \to t_0} \left( \left\Vert h_\alpha(t,\cdot)-h(t,\cdot) \right\Vert^{2}_{H^{-s}} - \left\Vert h_\alpha(t_0,\cdot)-h(t_0,\cdot) \right\Vert^{2}_{H^{-s}} \right)=0. \]

Once we have established this continuity property, there exists a time $0 < t_1 \leq T$ such that
	\[  \sup_{0\leq t \leq T} \left\Vert h_\alpha(t,\cdot)-h(t,\cdot) \right\Vert_{H^{-s}}=\left\Vert h_\alpha(t_1,\cdot)-h(t_1,\cdot) \right\Vert_{H^{-s}}.\]

Now, we prove the estimate (\ref{Above}). We write
\begin{equation}\label{Id}
	\left\Vert h_\alpha(t_1,\cdot)-h(t_1,\cdot) \right\Vert^{2}_{H^{-s}} = \int_{\Rt} \vert  e^{-\vert \xi \vert^{\alpha} t_1} -e^{-\vert \xi \vert^2 t_1}\vert^2 \frac{d \xi}{(1+\vert \xi \vert^2)^{s}}.
\end{equation} 

For $\xi \in \Rn \setminus \{0\}$ fixed, and for $1<\alpha < 2 +\delta$ (with $\delta>0$), we define the function  
\begin{equation}\label{function}
	f_\xi (\alpha)= e^{- t_1 \vert \xi \vert^\alpha},
\end{equation}
and by computing its derivative with respect to the variable $\alpha$, we get 
\begin{equation}\label{function-der}
f^{'}_{\xi}(\alpha)=- t_1\,  e^{- t_1 \vert \xi \vert^\alpha} \, \vert \xi \vert^\alpha \ln (\vert \xi \vert)
\end{equation}
Then, by the mean value theorem (in the variable $\alpha$), we can write 
\[ \vert f_\xi(\alpha)-f_\xi(2) \vert \leq \Vert f^{'}_{\xi} \Vert_{L^{\infty}([1,2+\delta])}\, \vert 2-\alpha \vert.\] 

Furthermore, we have the following uniform estimate with respect to the variable $\xi$:
\begin{equation}\label{Estim-Tech}
\sup_{\xi \in \Rt}  \Vert f^{'}_{\xi} \Vert_{L^{\infty}([1,2+\delta])}   \leq C\, T.
\end{equation}

Indeed,  we write 
\begin{equation*}
\left\Vert  \Vert f^{'}_{\xi} \Vert_{L^{\infty}([1,2+\delta])} \right\Vert_{L^{\infty}(\Rt)}\leq  \left\Vert  \Vert f^{'}_{\xi} \Vert_{L^{\infty}([1,2+\delta])} \right\Vert_{L^{\infty}(\vert \xi \vert \leq 1)}+ \left\Vert  \Vert f^{'}_{\xi} \Vert_{L^{\infty}([1,2+\delta])} \right\Vert_{L^{\infty}(\vert \xi \vert >1)} =A+B,
\end{equation*}
where we estimate the terms $A$ and $B$, separately. For the term $A$, as we have $\vert \xi \vert \leq 1$, $1<\alpha<2+\delta$, and moreover, since $\ds{\lim_{\vert \xi \vert\to 0^{+}} \vert \xi \vert \ln(\vert \xi \vert)=0}$, we deduce the following control:
\[ A \leq  T \, \left( \sup_{\xi \in \R^3} e^{-t_1 \vert \xi \vert^{2+\delta}} \vert \xi \vert \ln (\vert \xi \vert)\right)\leq C\, T.\]

For the term $B$, since $\vert \xi \vert >1$, we obtain
\[B \leq T\, \left( \sup_{\xi \in \R^3} e^{-t_1 \vert \xi \vert} \vert \xi \vert^{2+\delta} \ln (\vert \xi \vert)\right)\leq C\, T. \]	 

Once we have the estimate (\ref{Estim-Tech}), we can write
\[  \vert f_\xi(\alpha) - f_\xi(2) \vert \leq C\, T (2-\alpha).\]

Finally, we get back to the identity (\ref{Id}) to get
\begin{equation*}
	\begin{split}
		\Vert h_\alpha(t_1, \cdot) - h(t_1,\cdot) \Vert^{2}_{H^{-s}} = &\,  \int_{\Rn} \vert f_\xi (\alpha) - f_\xi(2) \vert^2 \frac{d \xi}{(1+\vert \xi \vert^2)^s} \\
		\leq &\,  C\, T^2\, (2-\alpha)^2\,  \int_{\Rt}  \frac{d \xi}{(1+\vert \xi \vert^2)^s} \leq \, C_s \, T^2\, ( 2-\alpha )^2. 
	\end{split}
\end{equation*}

We prove now the estimate (\ref{Below}). We write 
\[ \sup_{0\leq t \leq T} \left\Vert h_\alpha(t,\cdot)-h(t,\cdot) \right\Vert^{2}_{H^{-s}} \geq  \int_{2<|\xi |<4}    \left\vert  e^{-\vert \xi \vert^{\alpha} \frac{T}{2}} -e^{-\vert \xi \vert^2 \frac{T}{2}}\right\vert^2 \frac{d \xi}{(1+\vert \xi \vert^2)^{s}},\]
where we must study the expression at the right-hand side. By  the function $f_{\xi}(\alpha)$ defined in (\ref{function}) (with $\frac{T}{2}$ instead of $t_1$) we have 
\[ \int_{2<|\xi |<4}    \left\vert  e^{-\vert \xi \vert^{\alpha} \frac{T}{2}} -e^{-\vert \xi \vert^2 \frac{T}{2}}\right\vert^2 \frac{d \xi}{(1+\vert \xi \vert^2)^{s}}= \int_{2<|\xi |<4}    \vert  f_\xi(\alpha) - f_\xi(2)\vert^2 \frac{d \xi}{(1+\vert \xi \vert^2)^{s}}.\]
Moreover, by the expression (\ref{function-der}) (always with $\frac{T}{2}$ instead of $t_1$) and the Taylor formula we get 
\begin{equation*}
\begin{split}
| f_\xi(\alpha) - f_\xi(2)| =&\, | f'_{\xi}(2) (2-\alpha) + o (2-\alpha)|\\
\geq &\, | f'_{\xi}(2) (2-\alpha)  - (- o (2-\alpha))|\\
\geq &\, |f'_{\xi}(2)| (2-\alpha) - |o (2-\alpha)|.
\end{split}
\end{equation*}
Then, for $\frac{|f'_{\xi}(2)|}{2}>0$ there exists $0<\varepsilon_1\ll 1$ such that for $1+\varepsilon_1<\alpha <2$ we have $|o (2-\alpha)| \leq \frac{|f'_{\xi}(2)|}{2}(2-\alpha)$, hence, we get back to the previous estimate to obtain 
\[ | f_\xi(\alpha) - f_\xi(2)| \geq  \frac{|f'_{\xi}(2)|}{2}(2-\alpha). \]
Once we have this estimate at our disposal, we get back to the last integral to write 
\[ \int_{2<|\xi |<4}    \vert  f_\xi(\alpha) - f_\xi(2)\vert^2 \frac{d \xi}{(1+\vert \xi \vert^2)^{s}} \geq \left(\frac{T}{2}\right)^2  \left( \min_{2<|\xi|<4}  \frac{|f'_{\xi}(2)|}{2} \right)^2 \int_{2<|\xi|<2} \frac{d\xi}{(1+|\xi|^2)^s} =c^{2}_{s} \left(\frac{T}{2}\right)^2  (2-\alpha)^2,\]
which yields the wished estimate (\ref{Below}).  Proposition \ref{Key-Lemma} is proven. \finpv
\section{Proof of Proposition \ref{Prop:LWP}}\label{Sec:LWP}
The proof is rather standard, so we will only detail the main estimates. For a time $0<T<+\infty$, we consider the Banach space $\mathcal{C}\big([0,T], H^s(\Rt)\big)$, endowed with its natural norm $\| \cdot \|_{L^{\infty}_{t} H^{s}_{x}}$. On the right-hand side of equation (\ref{Mild}), the linear term is straightforward to estimate, and we have $\ds{\left\| e^{-(-\Delta)^{\alpha/2} t} \,  \vu_{0,\alpha} \right\|_{L^{\infty}_{t} H^{s}_{x}} \leq \| \vu_{0,\alpha} \|_{H^s}}$.

Thereafter, for $0<t<T$ fixed, the bilinear term is estimated as follows 
\begin{equation*}
	\left\|  \int_{0}^{t} e^{-(-\Delta)^{\alpha/2} (t-\tau)} \, \P\left((\vu_\alpha \cdot \vec{\nabla})\vu_\alpha\right)(\tau, \cdot) d \tau \right\|_{H^s} \leq C \int_{0}^{t} \| \vec{\nabla} h_\alpha(t-\tau,\cdot) \|_{L^1} \, \| \vu_\alpha \otimes \vu_\alpha (\tau, \cdot) \|_{H^s} d \tau. 
\end{equation*}

From \cite[Lemma $2.2$]{Yu-Zhai}, we have $\ds{ \| \vec{\nabla} h_\alpha(t-\tau,\cdot) \|_{L^1} \leq C (t-\tau)^{-\frac{1}{\alpha}}}$. On the other hand, since $s>3/2$, using the product laws in Sobolev spaces, we can write $\ds{\| \vu_\alpha \otimes \vu_\alpha (\tau, \cdot) \|_{H^s}  \leq C \| \vu_\alpha (\tau,\cdot) \|^{2}_{H^s}}$.  We thus obtain
\begin{equation*}
	C \int_{0}^{t} \| \vec{\nabla} h_\alpha(t-\tau,\cdot) \|_{L^1} \, \| \vu_\alpha \otimes \vu_\alpha (\tau, \cdot) \|_{H^s} d \tau \leq C \left( \int_{0}^{t} (t-\tau)^{-\frac{1}{\alpha}}\right) \| \vu_\alpha \|^{2}_{L^{\infty}_{t} H^{s}_{x}} \leq C\, \frac{T^{1-\frac{1}{\alpha}}}{1-\frac{1}{\alpha}}\,  \| \vu_\alpha \|^{2}_{L^{\infty}_{t} H^{s}_{x}}.
\end{equation*}

The existence and uniqueness of a mild solution $\vu_\alpha$ follows from Picard's iterative schema, as long as $\ds{4 C \| \vu_{0,\alpha} \|_{H^s} \frac{T^{1-\frac{1}{\alpha}}}{1-\frac{1}{\alpha}} <1}$,  which defines the time $T_\alpha$ as in (\ref{Talpha}). Proposition \ref{Prop:LWP} is proven. \finpv 

\section{From non-local to local Navier-Stokes equations} \label{Sec:NL-Ll}
In the following, $C>0$ denotes a generic constant that may change in each line, but it does not depend on the parameter $\alpha$. 
\subsection{Proof of Theorem \ref{Th-Main}}
For a time $0<T \leq T_0$ fixed,  we write 
\begin{equation}\label{estim-base}
	\begin{split}
		&\, \sup_{0 \leq t \leq T} \Vert \vu_{\alpha}(t,\cdot)-\vu_2 (t,\cdot)\Vert_{L^{\infty}} \\
		\leq &\,  \sup_{0\leq t \leq T} \left\Vert e^{-(-\Delta)^{\alpha/2} t} \vu_{0,\alpha} - e^{\Delta t} \vu_{0,2} \right\Vert_{L^{\infty}}  \\
		&\, +  \sup_{0\leq t \leq T}\left\Vert \int_{0}^{t} e^{-(-\Delta)^{\alpha/2} (t-\tau)} \P\left((\vu_\alpha \cdot \vec{\nabla})\vu_\alpha\right)(\tau, \cdot) d \tau - \int_{0}^{t} e^{\Delta (t-\tau)} \P\left((\vu_2 \cdot \vec{\nabla})\vu_2 \right)(\tau, \cdot) d \tau \right\Vert_{L^{\infty}}\\
		=& \, I_{\alpha}+J_{\alpha}.
	\end{split} 
\end{equation}

We begin by estimating each term on the right.  For the term $I_\alpha$, we get 
\begin{equation}\label{estim-lin-01}
	\begin{split}
		I_\alpha \leq  &\, \sup_{0\leq t \leq T} \left\Vert \left( e^{-(-\Delta)^{\alpha/2} t} - e^{\Delta t} \right) \vu_{0,\alpha} \right\Vert_{L^\infty} +  \sup_{0\leq t \leq T} \left\Vert e^{\Delta t} \left( \vu_{0,\alpha} - \vu_{0,2}\right) \right\Vert_{L^\infty} \\
		= &\, \sup_{0\leq t \leq T} \left\Vert \Big(  h_\alpha(t,\cdot)-h(t,\cdot) \Big) \ast \vu_{0,\alpha} \right\Vert_{L^\infty} +  \sup_{0\leq t \leq T} \left\Vert h(t,\cdot) \ast \left( \vu_{0,\alpha} - \vu_{0,2}\right) \right\Vert_{L^\infty} \\
		= &\, I_{\alpha,1}+I_{\alpha,2}.  
	\end{split}
\end{equation}

Afterwards, to estimate the term $I_{\alpha,1}$, one can apply the Bessel potential operators $(1-\Delta)^{-s/2}$ and $(1-\Delta)^{s/2}$ to deduce  
\begin{equation*}
	I_{\alpha,1} = \sup_{0\leq t \leq T}\left\Vert (1-\Delta)^{-s/2}\Big( h_\alpha(t,\cdot)-h(t,\cdot) \Big) \ast (1-\Delta)^{s/2} \vu_{0,\alpha} \right\Vert_{L^{\infty}}.
\end{equation*}

Thus, thanks to Young's inequalities (with $1+1/\infty=1/2+1/2$), we can write
\begin{equation}\label{estim01}
	\begin{split}
		I_{\alpha,1} &\leq  C\,  \sup_{0\leq t \leq T} \left(\left\Vert (1-\Delta)^{-s/2} \Big( h_\alpha(t,\cdot)-h(t,\cdot) \Big) \right\Vert_{L^{2}}\, \left\Vert (1-\Delta)^{s/2} \vu_{0,\alpha} \right\Vert_{L^{2}}\right) \\
		& \leq  C\left(  \sup_{0\leq t \leq T} \left\Vert  h_\alpha(t,\cdot)-h(t,\cdot) \right\Vert_{H^{-s}} \right)\, \left(  \sup_{1+\varepsilon < \alpha < 2 } \Vert \vu_{0,\alpha} \Vert_{H^{s}} \right),
	\end{split}
\end{equation}
where each of the terms above must be estimated separately. Note that, for the first term on the right-hand side, it is natural to apply the estimate (\ref{Above}) proven in  Proposition \ref{Key-Lemma}, whereas the second term on the right-hand side can be controlled by the fact that the family $\ds{(u_{0,\alpha})_{1+\varepsilon<\alpha <2}}$ is bounded in $H^s(\mathbb{R}^3)$.%(\ref{Conv-Data})  

\medskip 

Therefore,  the term $I_{\alpha,1}$ given in  (\ref{estim-lin-01}) can be estimated as follows 
\begin{equation}\label{Lim1-01}
	I_{\alpha,1} \leq C\, T\, (2-\alpha).
\end{equation}

It is now time to study the term $I_{\alpha,2}$ in (\ref{estim-lin-01}).  By Young's inequalities (with $1+ 1/\infty= 1 + 1/ \infty)$, the well-known properties of the heat kernel, and the assumption in (\ref{Conv-Rate-Data}), we have
\begin{equation}\label{Lim-02}
	I_{\alpha,2} \leq {\bf c} (2-\alpha)^\kappa.
\end{equation}

Consequently, we set the constant ${\bf C}_1=\max(C, {\bf c})$, and by using equations (\ref{Lim1-01}) and (\ref{Lim-02}), we can derive the following estimate
\begin{equation}\label{Lin}
	I_\alpha \leq {\bf C}_1\,(1+T)\,\Big( (2+\alpha)+ (2-\alpha)^\kappa\Big).
\end{equation}

\medskip

Similarly, the term $J_{\alpha}$ in (\ref{estim-base}) can also be studied separately.
\begin{equation}\label{Estim-J}
	\begin{split}
		J_{\alpha} \leq &  \sup_{0\leq t \leq T}\left\Vert \int_{0}^{t} h_\alpha(t-\tau,\cdot)\ast \P\left((\vu_\alpha \cdot \vec{\nabla})\vu_\alpha\right)(\tau, \cdot) d \tau - \int_{0}^{t} h(t-\tau,\cdot) \ast \P\left((\vu_\alpha \cdot \vec{\nabla})\vu_\alpha \right)(\tau, \cdot) d \tau \right\Vert_{L^{\infty}}\\
		&+  \sup_{0\leq t \leq T}\left\Vert \int_{0}^{t} h(t-\tau,\cdot) \ast \P\left((\vu_\alpha \cdot \vec{\nabla})\vu_\alpha\right)(\tau, \cdot) d \tau - \int_{0}^{t} h(t-\tau,\cdot) \ast \P\left((\vu_2 \cdot \vec{\nabla})\vu_2 \right)(\tau, \cdot) d \tau \right\Vert_{L^{\infty}}\\
		\leq & \sup_{0\leq t \leq T}  \left\Vert \int_{0}^{t}\Big(h_\alpha(t-\tau,\cdot)-h(t-\tau,\cdot)\Big)\ast \P\left((\vu_\alpha \cdot \vec{\nabla})\vu_\alpha\right)(\tau, \cdot) d \tau\right\Vert_{L^{\infty}} \\
		&+ \sup_{0\leq t \leq T} \left\Vert \int_{0}^{t}  h(t-\tau,\cdot)\ast \P\left((\vu_\alpha \cdot \vec{\nabla})\vu_\alpha - (\vu_2 \cdot \vec{\nabla})\vu_2 \right)(\tau, \cdot) d \tau \right\Vert_{L^{\infty}}=J_{\alpha,1}+J_{\alpha,2}.
	\end{split}
\end{equation}
 
For the term $J_{\alpha,1}$, we can leverage the properties of Leray's projector $\P$, and once again, we apply the operators $(1-\Delta)^{-s/2}$ and $(1-\Delta)^{s/2}$, along with Young's inequalities (with $1+1/\infty= 1/2 + 1/2$), to get the following estimates 
\begin{equation}\label{estim-J1}  
	\begin{split}
		J_{\alpha,1}
		& \leq  \sup_{0\leq t \leq T} \left( \int_{0}^{t} \left\Vert \Big(h_\alpha(t-\tau,\cdot)-h(t-\tau,\cdot)\Big)\ast \P\left(\text{div}(\vu_\alpha \otimes \vu_\alpha)\right)(\tau, \cdot) \right\Vert_{L^{\infty}} d\tau \right) \\
		& \leq  \sup_{0\leq t \leq T} \left(\int_{0}^{t} \left\Vert \P\Big(\vec{\nabla}h_\alpha(t-\tau,\cdot) - \vec{\nabla}h(t-\tau,\cdot)\Big)\right\Vert_{H^{-s}}\, \left\Vert\left(\vu_\alpha \otimes \vu_\alpha\right)(\tau, \cdot) \right\Vert_{H^s} d\tau \right) \\
		&\leq \sup_{0\leq t \leq T} \left( \int_{0}^{t} \left\Vert \vec{\nabla}h_\alpha(t-\tau,\cdot) - \vec{\nabla}h(t-\tau,\cdot)\right\Vert_{H^{-s}}\, \left\Vert \left(\vu_\alpha \otimes \vu_\alpha\right)(\tau, \cdot) \right\Vert_{H^s} d\tau  \right) \\
		&\leq T \left( \sup_{0\leq t \leq T} \left\Vert \vec{\nabla}h_\alpha(t,\cdot) -  \vec{\nabla}h(t,\cdot)\right\Vert_{H^{-s}}\right) \left( \sup_{0\leq t \leq T} \Vert \left(\vu_\alpha \otimes \vu_\alpha\right)(t, \cdot) \Vert_{H^s}\right).
	\end{split}
\end{equation}

To control the first term on the right-hand side, we can adapt Lemma \ref{Key-Lemma} to the function $f_{\xi}(\alpha) = i\xi_j e^{-t|\xi|\alpha}$, with $j = 1, 2, 3$, this manner, we obtain
\begin{align*}
	\sup_{0\leq t \leq T} \left\Vert \vec{\nabla}h_\alpha(t,\cdot) -  \vec{\nabla}h(t,\cdot)\right\Vert_{H^{-s}} \leq C T ( 2 - \alpha).
\end{align*}

For the remaining term on the right-hand side, we prove that there exists a constant ${\bf C}_2={\bf C}_2(\vu_{2,0}, \varepsilon)$ that is sufficiently large and depends only on $\vu_{0,2}$ and $\varepsilon$, such that the following uniform estimate holds:
\begin{equation*}
	\sup_{1+\varepsilon < \alpha <2} \sup_{0\leq t \leq T} \| (\vu_{\alpha}\otimes \vu_{\alpha})(t,\cdot) \|_{H^s} \leq {\bf C}_2.
\end{equation*}

Indeed, recall that the solution $\vu_\alpha \in \mathcal{C}\big([0,T], H^s(\Rt)\big)$ obtained in Proposition \ref{Prop:LWP} by the Picard's iterative argument verifies
\[ \sup_{0\leq t \leq T} \Vert \vu_\alpha(t,\cdot) \Vert_{H^s} \leq \sup_{0\leq t \leq T_\alpha} \Vert \vu_\alpha(t,\cdot) \Vert_{H^s} \leq  C\,  \Vert \vu_{0,\alpha} \Vert_{H^s}, \quad \text{where } T\leq T_0 \leq T_\alpha. \]

Moreover, based on the assumption (\ref{Conv-Data}), we have $\ds{\sup_{1+\varepsilon < \alpha < 2} \| \vu_{0,\alpha} \|_{H^s}}\leq {\bf C}_2$. Then, we obtain 
\begin{equation}\label{Estim-unif-1}
	\sup_{1+\varepsilon < \alpha <2} \,  \sup_{0\leq t \leq T} \| \vu_{\alpha}(t,\cdot) \|_{H^s} \leq {\bf C}_2.
\end{equation}

Thus, the desired estimate follows from the fact that $s>3/2$ and, using the product laws in Sobolev spaces,  we can write 
\begin{equation*}
	\sup_{0 \leq t \leq T} \Vert \left(\vu_\alpha \otimes \vu_\alpha\right)(t,\cdot)\Vert_{H^s} \leq  C\,  \sup_{0 \leq t \leq T} \Vert \vu_\alpha(t,\cdot)\Vert^{2}_{H^s} 
	\leq \, C \left( \sup_{0 \leq t \leq T} \Vert \vu_\alpha(t,\cdot)\Vert_{H^s}  \right)^{2} \leq {\bf C}_2.
\end{equation*}

Returning to estimate (\ref{estim-J1}), the above inequality allows us to write
\begin{equation}\label{Lim2}
	J_{\alpha,1} \leq {\bf C}_2\, T^2 \, \vert 2-\alpha \vert \leq  {\bf C}_2\, T^2\, \Big((2-\alpha)+ (2-\alpha)^\kappa\Big).
\end{equation}

\medskip 

Subsequently, we study the term $J_{\alpha,2}$ given in (\ref{Estim-J}). For this propose we combine  Leray's projector $\P$ properties and Young's inequalities (with $1+ 1/\infty= 1 + 1/ \infty$) as follows
\begin{equation}\label{eq:30}
	\begin{split}
		J_{\alpha,2} &\leq \sup_{0\leq t \leq T} \int_{0}^{t}  \left\Vert h(t-\tau,\cdot)\ast \P\big(\text{div}(\vu_\alpha \otimes \vu_\alpha) - \text{div}(\vu_2 \otimes \vu_2)\big)(\tau, \cdot) \right\Vert_{L^{\infty}} d \tau\\
		&\leq C\, \sup_{0\leq t \leq T} \int_{0}^{t} \Vert \nabla h(t-\tau,\cdot)\Vert_{L^1} \Vert \P \big((\vu_\alpha \otimes \vu_\alpha) - (\vu_2 \otimes \vu_2)\big)(\tau,\cdot) \Vert_{L^\infty} d\tau.
	\end{split}
\end{equation}

Due to the well-known properties of the heat kernel $h(t,\cdot)$, we have $\Vert \nabla h(t-\tau,\cdot)\Vert_{L^1} \leq C (t-\tau)^{-1/2}$. Meanwhile, to estimate the term $\Vert \P \big((\vu_\alpha \otimes \vu_\alpha) - (\vu_2 \otimes \vu_2)\big)(\tau,\cdot) \Vert_{L^\infty}$, we make use of Leray's projector $\P$ properties, the uniform estimate  inequality (\ref{Estim-unif-1}) and the fact that $s>3/2$. Thus,
\begin{equation*}
	\begin{split}
		& \, 	\Vert \P \big((\vu_\alpha \otimes \vu_\alpha) - (\vu_2 \otimes \vu_2)\big)(\tau,\cdot) \Vert_{L^\infty} \\%= \left\Vert \P \big( (\vu_{\alpha}-\vu_{2}) \otimes (\vu_{\alpha}+\vu_{2}) \big) (\tau,\cdot) \right\Vert_{L^\infty} \\
		= & \, \left\Vert \big(\vu_{\alpha}(\tau,\cdot)-\vu_{2}(\tau,\cdot)\big) \otimes \P (\vu_{\alpha}+\vu_{2}) (\tau,\cdot) \right\Vert_{L^\infty} \\
		\leq &\, \left\Vert \vu_{\alpha}(\tau,\cdot)-\vu_{2}(\tau,\cdot)\right\Vert_{L^\infty} \big( \left\Vert \P (\vu_{\alpha})(\tau,\cdot)\right\Vert_{L^\infty}+ \left\Vert \P (\vu_{2})(\tau,\cdot) \right\Vert_{L^\infty}\big)  \\
		\leq & \, \left\Vert \vu_{\alpha}(\tau,\cdot)-\vu_{2}(\tau,\cdot)\right\Vert_{L^\infty} \big( \left\Vert\vu_{\alpha}(\tau,\cdot)\right\Vert_{H^s}+ \left\Vert\vu_{2}(\tau,\cdot) \right\Vert_{H^s}\big)  \\
		\leq &  \,  {\bf C}_2\,  \Vert \vu_{\alpha}(\tau,\cdot)-\vu_{2}(\tau,\cdot) \Vert_{L^\infty}.
	\end{split}
\end{equation*}

These last two estimations allow us to control \eqref{eq:30} as follows
\begin{equation}\label{Lim3}
	J_{\alpha,2} \leq  \, {\bf C}_2\, \sup_{0\leq t \leq T} \, \int_{0}^{t}(t-\tau)^{-1/2} \, \Vert  \vu_{\alpha}(\tau,\cdot)-\vu_{2}(\tau,\cdot) \Vert_{L^\infty} d\tau \leq \,  {\bf C}_2\, T^{1/2}\, \left( \sup_{0\leq t \leq T} \Vert  \vu_{\alpha}(t,\cdot)-\vu_{2}(t,\cdot) \Vert_{L^\infty}  \right).
\end{equation}

Thus far, we have controlled the terms $I_\alpha$, $J_{\alpha,1}$ and $J_{\alpha,2}$ in (\ref{Lin}), (\ref{Lim2}), and (\ref{Lim3}), respectively.  We set the constant ${\bf C}= \max({\bf C}_1, {\bf C}_2)$, and we get back to (\ref{estim-base}) to write
\begin{equation*}
	\begin{split}
		\sup_{0\leq t \leq T} \Vert \vu_\alpha(t,\cdot)- \vu_2(t,\cdot)\Vert_{L^\infty} \leq &\,  I_\alpha + J_{\alpha,1}+J_{\alpha, 2}  \leq I_\alpha + J_{\alpha,1} +  {\bf C}\,  T^{1/2}\, \left( \sup_{0\leq t \leq T} \Vert  \vu_{\alpha}(t,\cdot)-\vu_{2}(t,\cdot) \Vert_{L^\infty}  \right). 
	\end{split}
\end{equation*}

In the above estimate, we set a time $0<T_1\leq T$ such that $\ds{	{\bf C}\, T^{1/2}_{1} \leq \frac{1}{2}}$. This way, we derive the following control:
\[  \sup_{0\leq t \leq T_1} \Vert \vu_\alpha(t,\cdot)- \vu_2(t,\cdot)\Vert_{L^\infty}  \leq I_\alpha + J_{\alpha,1}+  \frac{1}{2}\, \left( \sup_{0\leq t \leq T_1} \Vert  \vu_{\alpha}(t,\cdot)-\vu_{2}(t,\cdot) \Vert_{L^\infty}  \right),\]
and we can write  
\[ \frac{1}{2} \sup_{0\leq t \leq T_1} \Vert \vu_\alpha(t,\cdot)- \vu_2(t,\cdot)\Vert_{L^\infty} \leq I_\alpha + J_{\alpha,1}.\]

Then, by (\ref{Lin}) and (\ref{Lim2}), we obtain
\begin{equation}\label{Iteration-time}
%	\begin{split}
		\sup_{0\leq t \leq T_1} \Vert \vu_\alpha(t,\cdot)- \vu_2(t,\cdot)\Vert_{L^\infty}  % & \, {\bf C}(1+T_1+T^{2}_{1}) \, \max\left( (2-\alpha)^\gamma, (2-\alpha)\right)\\
		\leq \,  {\bf C}(1+T_1+T^{2}_{1}) \, \Big( (2-\alpha)+(2-\alpha)^\kappa \Big).  
%	\end{split}
\end{equation}

By iterative application of this argument up to time $T_0 > 0$, we have
\begin{equation}\label{Key-estimate}
	\sup_{0\leq t \leq T_0} \Vert \vu_\alpha(t,\cdot)- \vu_2(t,\cdot)\Vert_{L^\infty} \leq  \, {\bf C}(1+T_0+T^{2}_{0}) \, \Big((2-\alpha)+ (2-\alpha)^\kappa\Big). 
\end{equation}

To finish the proof of Theorem \ref{Th-Main}, we shall prove that estimate (\ref{Key-estimate}) yields
\begin{equation}\label{Pressures}
	\sup_{0\leq t \leq T_0} \Vert p_\alpha(t,\cdot)- p_2(t,\cdot)\Vert_{BMO} \leq  \, {\bf C}(1+T_0+T^{2}_{0}) \, \Big((2-\alpha)+ (2-\alpha)^\kappa\Big).
\end{equation}

Indeed, using expression (\ref{Pressure}), the estimate $\| \mathcal{R}_i f \|_{BMO} \leq C \| f \|_{L^\infty}$ (see, for instance, \cite[Theorem $6.2$]{PGLR0}), and the uniform estimate (\ref{Estim-unif-1}), for $0<t\leq T_0$, we write 
\begin{equation*}
	\begin{split}
		\| p_\alpha(t,\cdot)-p_2(t,\cdot)\|_{BMO} \leq &\,  C \| \vu_{\alpha} \otimes \vu_{\alpha} (t,\cdot) - \vu_{2} \otimes \vu_{2} (t,\cdot) \|_{L^\infty} \\
		\leq &\,  C \| \vu_{\alpha}(t,\cdot) - \vu_2(t,\cdot) \|_{L^\infty} \big( \| \vu_{\alpha}(t,\cdot)\|_{L^\infty}+ \| \vu_{2}(t,\cdot)\|_{L^\infty} \big) \\
		\leq & \, {\bf C} \, \| \vu_{\alpha}(t,\cdot) - \vu_2(t,\cdot) \|_{L^\infty},
	\end{split}
\end{equation*}
which yields (\ref{Pressures}). Theorem \ref{Th-Main} is now proven. \finpv

\subsection{Proof of Corollary \ref{Corollary:global}} 
First, observe that in the case of global in time mild solutions (under the assumption (\ref{Small-data}))  we can iterate (\ref{Iteration-time}) to obtain (\ref{Conv-Rate-Sol-global}) for any time $0<T<+\infty$.

\medskip

On the other hand, by the uniform estimate (\ref{Estim-unif-1}), the limit solution $\vu_2$ verifies  $\vu \in L^{\infty}_{loc}([0,+\infty), H^s(\Rt))$, which yields that $\vu_2$  belongs to the energy space $L^{\infty}_{t}L^{2}_{x}\cap (L^{2}_{loc})_t \dot{H}^{1}_{x}$ and it also verifies an energy equality. Corollary \ref{Corollary:global} is proven. \finpv 

\subsection{Proof of Corollary \ref{Corollary:LpLq}} 
Remark that the family of initial data also belongs to the space $L^2(\Rt)$, and by well-known arguments, for $1+\varepsilon<\alpha <1$, we have 
\[ \| \vu_{\alpha}(t,\cdot) \|^{2}_{L^2} \leq  \| \vu_{0,\alpha} \|^{2}_{L^2} \leq C \| \vu_{0,\alpha}\|^{2}_{H^s}\leq {\bf C}_2.\] 

Estimate (\ref{Conv-Rate-Sol-Lp}) follows from a standard interpolation argument (in Lebesgue spaces) between the estimate above and (\ref{Conv-Rate-Sol}). Then, Corollary \ref{Corollary:LpLq} is proven. \finpv
%%%%%%%%%%%%%%%%%%%%%%%%%%%%%%%%%%%%%%%%%%%%%%%%%%%
\begin{appendices}
	
\section{Appendix}\label{AppendixA}
We now prove the lower bound (\ref{Lower-bound-time}). Using (\ref{Conv-Data}), we can set $0<\varepsilon \ll 1$ such that for all $1+\varepsilon < \alpha <2$, we have  $\ds{\left| \| \vu_{0,\alpha} \|_{H^s}- \| \vu_{0,2} \|_{H^s} \right| \leq \frac{1}{2} \| \vu_{0,2} \|_{H^s}}$. Thus, we obtain $\ds{\| \vu_{0,\alpha} \|_{H^s}  \leq \frac{3}{2}  \| \vu_{0,2} \|_{H^s}}$, and we can write 
\begin{equation*}
	\frac{1}{2} \left( \frac{1-\frac{1}{\alpha}}{4C \| \vu_{0,2} \|_{H^s}} \right)^{\frac{\alpha}{\alpha-1}} \leq T_\alpha, \quad 1+\varepsilon < \alpha < 2. 
\end{equation*}

Furthermore, the expression on the left-hand side is estimated from below by the quantity 
\begin{equation*}
	T_0= \frac{1}{2} \max\left[  \left( \frac{1-\frac{1}{1+\varepsilon}}{4C \| \vu_{0,2} \|_{H^s}} \right)^{\frac{2}{\varepsilon}}, \left( \frac{1-\frac{1}{1+\varepsilon}}{4C \| \vu_{0,2} \|_{H^s}} \right)^{1+\varepsilon}\right].  
\end{equation*}

Indeed, as we have $1+\varepsilon < \alpha < 2$, then we get $1-\frac{1}{1+\varepsilon}< 1-\frac{1}{\alpha}$, and write 
	\[ \frac{1}{2} \left( \frac{1-\frac{1}{1+\varepsilon}}{4 C \Vert \vu_{0,2} \Vert_{H^s}}\right)^{\frac{\alpha}{\alpha-1}} \leq \frac{1}{2}\left( \frac{1-\frac{1}{\alpha}}{4 C\Vert \vu_{0,2} \Vert_{H^s}} \right)^{\frac{\alpha}{\alpha-1}}.\]
	
Thereafter, for the sake of simplicity, we denote $\ds{A=\frac{1-\frac{1}{1+\varepsilon}}{4 C \Vert \vu_{0,2} \Vert_{H^s}}}$,  and we have 
\[ \frac{1}{2} A^{\frac{\alpha}{\alpha-1}} \leq \frac{1}{2}\left[ \frac{1-\frac{1}{\alpha}}{4C \Vert \vu_{0,2} \Vert_{H^s}} \right]^{\frac{\alpha}{\alpha-1}}. \]

Now, let us study the expression $\frac{\alpha}{\alpha-1}$.  Since $1+\varepsilon < \alpha <2$, then we get $1+\varepsilon < \frac{\alpha}{\alpha-1} < \frac{2}{\varepsilon}$.  Thus, on one hand, if the quantity $A$ above verifies $A<1$, then we have $\ds{\frac{1}{2}A^{\frac{2}{\varepsilon}} \leq  \frac{1}{2}A^{\frac{\alpha}{\alpha-1}}}$.  On the other hand, if the quantity $A$ satisfies $1 \leq A$, then we obtain $\ds{\frac{1}{2}A^{1+\varepsilon} \leq  \frac{1}{2}A^{\frac{\alpha}{\alpha-1}}}$.

\section{Appendix}\label{AppendixB}
For $1<\alpha,\beta \leq 2$, we consider the initial value problem for the MHD equations
	\begin{equation}\label{MHD} 
		\begin{cases}\vspace{2mm}
			\partial_t \vu = -(-\Delta)^{\alpha/2} \vu - (\vu \cdot \vec{\nabla})\vu + (\vb \cdot \vec{\nabla} )\vb -\vec{\nabla} p, \quad \text{div}(\vu)=0, \\ \vspace{2mm}
			\partial_t \vb = -(-\Delta)^{\beta/2} \vb - (\vu \cdot \vec{\nabla})\vb + (\vb \cdot \vec{\nabla} )\vu, \quad \text{div}(\vb)=0, \\ 
			\vu_0(t,\cdot)=\vu_0, \quad \vb(0,\cdot)=\vb_0,
			\end{cases}
	\end{equation}
	where $\vu: [0,+\infty) \times \Rt \to \Rt$ and $p:\vb: [0,+\infty) \times \Rt \to \R $ always denote the velocity and the pressure of the fluid,   $\vb: [0,+\infty) \times \Rt \to \Rt $ is the magnetic field, and  $\vu_0, \vb_0 : \Rt \to \Rt$ are the divergence-free initial data. 
	
\medskip

Recall that mild solutions to  the system  (\ref{MHD}) write down as 
\begin{equation*}\label{Mild1}
\begin{split}
\vu(t,\cdot)=e^{-(-\Delta)^{\alpha/2} t} \, \vu_{0} - \underbrace{\int_{0}^{t} e^{-(-\Delta)^{\alpha/2} (t-\tau)} \, \P\left((\vu \cdot \vec{\nabla})\vu  \right)(\tau, \cdot) d \tau}_{B_1(\vu,\vu)}  + \underbrace{ \int_{0}^{t} e^{-(-\Delta)^{\alpha/2} (t-\tau)} \, \P\left( (\vb \cdot \vec{\nabla})\vb \right)(\tau, \cdot) d \tau}_{B_{2}(\vb,\vb)}, 
\end{split}
\end{equation*}
\begin{equation*}\label{Mild2}
\begin{split}
\vb(t,\cdot)= e^{-(-\Delta)^{\beta/2} t} \, \vb_{0} - \underbrace{\int_{0}^{t} e^{-(-\Delta)^{\beta/2} (t-\tau)} \, \P\left((\vu \cdot \vec{\nabla})\vb  \right)(\tau, \cdot) d \tau}_{B_3(\vu,\vb)}+ \underbrace{\int_{0}^{t} e^{-(-\Delta)^{\beta/2} (t-\tau)} \, \P\left((\vb \cdot \vec{\nabla})\vu \right)(\tau, \cdot) d \tau}_{B_4(\vb,\vu)},
\end{split}
\end{equation*}
and the pressure term is related to the velocity $\vu$ and the magnetic field $\vb$ by the well-known expression 
\begin{equation*}\label{PressureMHD}
p=\sum_{i,j=1}^{3} \mathcal{R}_{i}\mathcal{R}_{j} (u_{i} \, u_{j}+b_{i} \, b_{j}).
\end{equation*}

\medskip

As in Proposition \ref{Prop:LWP}, the existence of local-in-time mild $H^{s}$-solutions  (with $s>3/2$) is rather classical, and we can state the following result adapted to the coupled system (\ref{MHD}):
\begin{Proposition}\label{Prop:LWP-MHD} Let $1<\alpha,\beta \leq 2$ be fixed. Let $s>3/2$ and let $\vu_{0,\alpha}, \vb_{0,\beta} \in H^s(\Rt)$ be a divergence-free initial datum. There exists a time 
		\begin{equation}\label{Talpha-beta}
		0<T_{\alpha,\beta} = \frac{1}{2} \min\left[\left( \frac{1- \frac{1}{\alpha}}{4C \| \vu_{0,\alpha} \|_{H^s}} \right)^{\frac{\alpha}{\alpha-1}},\left( \frac{1- \frac{1}{\beta}}{4C \| \vb_{0,\beta} \|_{H^s}} \right)^{\frac{\beta}{\alpha-1}}\right],
		\end{equation}	
		where $C>0$ is a generic constant, and there exists a unique mild solution $\vu_{\alpha,\beta}, \vb_{\alpha,\beta}$ to the system (\ref{MHD}), such that 
		\[ \vu_{\alpha,\beta}, \ \vb_{\alpha,\beta} \in \mathcal{C}\big([0,T_{\alpha,\beta}], H^s(\Rt)\big) \quad \mbox{and} \quad  p_{\alpha,\beta}  \in \mathcal{C}\big([0,T_{\alpha,\beta}], H^s(\Rt)\big).\]
	\end{Proposition}	

As before, we shall assume the following strong  convergence in the space $H^s(\Rt)$:
\begin{equation}\label{Conv-Data-MHD}
\vu_{0,\alpha} \to \vu_{0,2}  \quad \alpha \to 2, \quad   \vb_{0,\beta} \to \vb_{0,2} \quad \beta \to 2,
\end{equation}
which yields the existence of a parameter $0<\varepsilon\ll 1$ and a time $T_0$ (depending on $\varepsilon$), such that the entire family of solutions $(\vu_{\alpha,\beta}, \vb_{\alpha,\beta})_{1+\varepsilon<\alpha,\beta\leq 2}$ is at least well-defined in the interval $[0,T_0]$.  See always Appendix \ref{AppendixA} for a rigorous explanation. 

\medskip

Then, the non-local to local convergence adapted to the MHD system reads as follows:
\begin{Theoreme}\label{Th-Main-MHD} Let $(\vu_{0,\alpha}, \vb_{0,\beta})_{1+\varepsilon<\alpha, \beta \leq 2}$ be an initial data family, where $\vu_{0,\alpha}, \vb_{0,\beta} \in H^s(\Rt)$ with $s>3/2$. Let $\ds{(\vu_{\alpha,\beta}, \vb_{\alpha,\beta},  p_{\alpha,\beta})_{1+\varepsilon<\alpha,\beta\leq 2} \subset  \mathcal{C}\big([0,T_0], H^s(\Rt)\big)}$ be the corresponding family of solutions to the system (\ref{MHD}), given by Proposition \ref{Prop:LWP-MHD}. 
	
	\medskip 
	
We assume the convergence given in (\ref{Conv-Data-MHD}), and we assume the estimates 
	\begin{equation}\label{Conv-Rate-Data-MHD}
	\|  \vu_{0,\alpha} - \vu_{0,2} \|_{L^\infty} \leq {\bf c_1} (2-\alpha)^{\kappa_1}, \quad  	\|  \vb_{0,\alpha} - \vb_{0,2} \|_{L^\infty} \leq {\bf c_2} (2-\beta)^{\kappa_2}, \quad \mbox{with} \ \ 0<\kappa_1,\kappa_2,
	\end{equation}
and where ${\bf c_1}, {\bf c_2}>0$ are generic constants.  

Then, there exists a constant $0<{\bf C}(T_0)\sim 1+T_0 +T^{2}_{0}$, depending on $\vu_{0,2}, \varepsilon$, ${\bf c}_1$,  ${\bf c}_2$, and the time $T_0$, such that for all $1+\varepsilon < \alpha < 2$ the following estimate holds:
	\begin{equation}\label{Conv-Rate-Sol-MHD}
	\begin{split}
	\sup_{0\leq t \leq T_0} & \big( \| \vu_\alpha(t,\cdot) - \vu_2(t,\cdot)\|_{L^\infty} + \| \vb_\alpha(t,\cdot) - \vb_2(t,\cdot)\|_{L^\infty} + \| p_\alpha(t,\cdot) - p_2(t,\cdot)\|_{BMO}\big) \\
	 \leq &  {\bf C}(T_0)  \max\Big(  (2-\alpha)+(2-\alpha)^{\kappa_1},  (2-\beta)+(2-\beta)^{\kappa_2}\Big). 
	\end{split}
	\end{equation}
\end{Theoreme}

As in Theorem \ref{Th-Main}, we observe that the convergence rate assumed for initial data in (\ref{Conv-Rate-Data-MHD}) does not always propagate to solutions due to the prescribed convergence rate of the kernels $h_\alpha(t,\cdot)\to h(t,\cdot)$ and $h_\beta(t,\cdot)\to h(t,\cdot)$, when $\alpha, \beta \to 2^{-}$. 

\medskip

\pv We essentially follow the same lines in the proof of Theorem \ref{Th-Main}, so it is enough to provide a brief sketch. For a time $0<T\leq T_0$, we start by writing
\begin{equation*}
\begin{split}
&\, \sup_{0 \leq t \leq T} \left(  \Vert \vu_{\alpha}(t,\cdot)-\vu_2 (t,\cdot)\Vert_{L^{\infty}} +  \Vert \vb_{\beta}(t,\cdot)-\vb_2 (t,\cdot)\Vert_{L^{\infty}} \right)\\ 
\leq &\,  \sup_{0\leq t \leq T} \left( \left\Vert e^{-(-\Delta)^{\alpha/2} t} \vu_{0,\alpha} - e^{\Delta t} \vu_{0,2} \right\Vert_{L^{\infty}} + \left\Vert e^{-(-\Delta)^{\beta/2} t} \vu_{0,\beta} - e^{\Delta t} \vu_{0,2} \right\Vert_{L^{\infty}} \right) \\
&\,  +  \sup_{0\leq t \leq T}\left\Vert B_1(\vu_{\alpha,\beta}, \vu_{\alpha,\beta}) - B_1(\vu_2,\vu_2) \right\Vert_{L^{\infty}} +  \sup_{0\leq t \leq T}\left\Vert B_2(\vb_{\alpha,\beta}, \vb_{\alpha,\beta}) - B_2(\vb_2,\vb_2) \right\Vert_{L^{\infty}} \\
&\,  +  \sup_{0\leq t \leq T}\left\Vert B_3(\vu_{\alpha,\beta}, \vb_{\alpha,\beta}) - B_3(\vu_2,\vb_2) \right\Vert_{L^{\infty}} +  \sup_{0\leq t \leq T}\left\Vert B_4(\vb_{\alpha,\beta}, \vu_{\alpha,\beta}) - B_4(\vb_2,\vu_2) \right\Vert_{L^{\infty}}  \\
=& \, I_{\alpha}+ I_{\beta} + \sum_{i=1}^{4} J_{\alpha,\beta,i}.
\end{split}  
\end{equation*}

Terms $I_{\alpha}$ and $I_{\beta}$ are estimated as in (\ref{Lin}), and we have
\begin{equation*}
I_{\alpha}+ I_{\beta} \leq {\bf C}_1 (1+T) \max\Big(  (2-\alpha)+(2-\alpha)^{\kappa_1},  (2-\beta)+(2-\beta)^{\kappa_2}\Big).
\end{equation*}

Thereafter, for $i=1,\cdots,4$, the  terms $J_{\alpha,\beta,i}$ are estimated as in (\ref{Lim2}) and (\ref{Lim3}) to obtain
\begin{equation*}
\begin{split}
\sum_{i=1}^{4} J_{\alpha,\beta,i} \leq  &\, {\bf C}_2 T^2  \max\Big(  (2-\alpha)+(2-\alpha)^{\kappa_1},  (2-\beta)+(2-\beta)^{\kappa_2}\Big) \\
&\,  + {\bf C}_2 T^{1/2}  \sup_{0 \leq t \leq T} \left(  \Vert \vu_{\alpha}(t,\cdot)-\vu_2 (t,\cdot)\Vert_{L^{\infty}} +  \Vert \vb_{\beta}(t,\cdot)-\vb_2 (t,\cdot)\Vert_{L^{\infty}} \right), 
\end{split}
\end{equation*}
and we conclude with the proof as in the end of the proof of Theorem \ref{Th-Main}.  \finpv
\end{appendices}

%%%%%%%%%%%%%%%%%%%%%%%%%%%%%%%%%%%%%%%%%%%%%%%%%%%

\section*{Statements and Declarations}
Data sharing does not apply to this article as no datasets were generated or analyzed during the current study. This work has not received any financial support. In addition, the authors declare that they have no conflicts of interest, and all of them have equally contributed to this paper.

\end{document}